\documentclass[10pt]{article}

\usepackage[draft]{graphics}
\usepackage{graphicx}          	 
\usepackage{bm}                 
\usepackage{amsmath}             
\usepackage{amssymb}
\usepackage{amsfonts}             
\usepackage{verbatim}           
\usepackage{amsthm}             

\usepackage{mathtools}
\usepackage{relsize}

\usepackage[colorlinks,citecolor=blue]{hyperref}

\usepackage{makeidx}

\theoremstyle{plain}             

\newtheorem{theorem}{Theorem}[section]

\theoremstyle{definition}
\newtheorem{example}[theorem]{Example}

\newtheorem{remark}[theorem]{Remark}

\makeatletter
\@addtoreset{equation}{chapter}

\makeatother

\def\protectbold#1{\protect{\boldmath{$#1$}}}

\def\eqref#1{(\ref{#1})}
\def\dsp{\displaystyle}
\def\Frac#1#2{\frac
{
 {\raise.6ex
 \hbox{$\displaystyle#1$}}
}
{
 {\lower.6ex
 \hbox{$\displaystyle#2$}}
 }
}

\numberwithin{equation}{section}

\newcommand{\arcsinh}{{\rm arcsinh}}
\newcommand{\arccosh}{{\rm arccosh}}

\def\bigOxe{\sqcup \kern-2.3mm \sqcap}

\def\eoexample{{\unskip\nobreak\hfil\penalty50	
\hskip2em\hbox{}\nobreak\hfil$\diamondsuit$
\parfillskip=0pt\finalhyphendemerits=0\medbreak}}

\def\eoremark{{\unskip\nobreak\hfil\penalty50	
\hskip2em\hbox{}\nobreak\hfil$\triangle$
\parfillskip=0pt\finalhyphendemerits=0\medbreak}}

\makeatother

\def\dsp{\displaystyle}
\def\Frac#1#2{\frac
{
 {\raise.6ex
 \hbox{$\displaystyle#1$}}
}
{
 {\lower.6ex
 \hbox{$\displaystyle#2$}}
 }
}

\def\sign{{\rm sign}}

\input epsf

\def\CHFs#1#2#3{
{}_1F_1\left({a};{c};{z}\right)
}

\def\intr{\int_{-\infty}^\infty}
\def\bigO{{\cal O}}
\def\calC{{{\cal C}}}

\def\calL{{{\cal L}}}

\def\wt{\widetilde}

\def\CC{\mathbb C}
\def\RR{\mathbb R}

\def\Z z_+{{\mathbb Z^+}}      

\def\sign{{\rm sign}}

\def\Ai{{{\rm Ai}}}

\def\Gi{{{\rm Gi}}}
\def\Hi{{{\rm Hi}}}

\def\intr{\int_{-\infty}^\infty}

\def\calL{{\cal L}}

\def\phase{{\rm ph}}
\def\wt{\widetilde}
\def\tfrac#1#2{{{\lower.6ex
\hbox{$\scriptstyle#1$}}\over 
{\raise.7ex
\hbox{$\scriptstyle#2$}}}}

\def\arcsinh{{\rm arcsinh}}
\def\arccosh{{\rm arccosh}}
\def\arctanh{{\rm arctanh}}

\begin{document}
 \title{Numerical evaluation of Airy-type integrals arising in uniform asymptotic analysis}

\author{
A. Gil\\
Departamento de Matem\'atica Aplicada y CC. de la Computaci\'on.\\
ETSI Caminos. Universidad de Cantabria. 39005-Santander, Spain.\\
 \and
J. Segura\\
        Departamento de Matem\'aticas, Estadistica y 
        Computaci\'on,\\
        Universidad de Cantabria, 39005 Santander, Spain.\\
\and
N. M. Temme\\
IAA, 1825 BD 25, Alkmaar, The Netherlands.\footnote{Former address: Centrum Wiskunde \& Informatica (CWI), Science Park 123, 1098 XG Amsterdam,  The Netherlands}\\
}


\maketitle
\begin{abstract}
We describe a method to evaluate integrals that arise in the asymptotic analysis when two saddle points may be close together. These integrals,
which appear in problems from optics, acoustics or quantum mechanics as well as in a wide class of special functions, can be transformed into Airy-type integrals and we use the trapezoidal rule to compute these integrals numerically. The quadrature method, which remains valid when two saddle points coalesce, is illustrated with numerical examples.  
\end{abstract}

\medskip


\section{Introduction}\label{sec:intro}
We consider the numerical evaluation of  integrals of the form
\begin{equation}\label{eq:intro01}
F(\eta)=\frac{1}{2\pi i}\int_{\calC}e^{\frac13t^3-\eta t}f(t)\,dt \quad  \eta\in\RR.
\end{equation}
The contour $\calC$ runs from $\infty e^{-\frac13\pi i}$ to $\infty e^{+\frac13\pi i}$. The function $\phi(t)=\frac13t^3-\eta t$, with derivative $\phi^\prime(t)=t^2-\eta$, has two saddle points $\pm\sqrt\eta$, the zeros of $\phi^\prime(t)$. The function $f(t)$ is assumed to be analytic in a neighborhood of the contour~$\calC$.  

When we want to compute this integral, it will be convenient to choose the contour suitably, which means that we like to avoid contours where the dominant part of the integrand $e^{\phi(t)}$ is strongly oscillating.  We  already prescribed that the contour terminates in two valleys of the function $e^{\frac13t^3}$, but we can do more by trying to choose a path on which $\Im\phi(t)$ is constant, always verifying if the function $f(t)$ is allowing such a choice.

\begin{figure}
\begin{center}
\epsfxsize=10cm \epsfbox{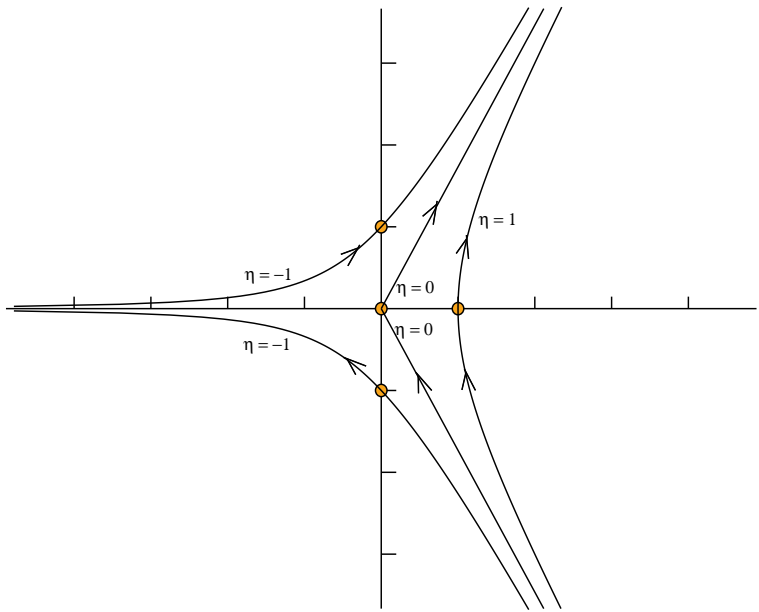}
\caption{
\label{fig:fig01} Saddle point contours $\calC$ for the integral in \eqref{eq:intro01} for $\eta=1$, $\eta=0$, and $\eta=-1$.}
\end{center}
\end{figure}

According to the classical methods of asymptotic analysis for contour integrals (see, for example,  \cite[Chapter 4]{Temme:2015:AMI}), the optimal choice is taking the path $\calC$ trough a saddle point, say, $\sqrt\eta$,  on which  $\Im\phi(t)=\Im\phi\left(\sqrt\eta\right)$. When $\eta>0$ this is indeed the best choice. In that case $\Im\phi\left(\sqrt\eta\right)=0$, and by writing $t=u+iv$, the equation that defines this optimal path $\calC$ is described by $\Im\phi(t)=v\left(u^2-\frac13v^2-\eta\right)=0$. The path that terminates in the valleys as assumed for the integral in \eqref{eq:intro01} is governed by 
\begin{equation}\label{eq:intro02}
u=\sqrt{\tfrac13v^2+\eta}, \quad v\in\RR.
\end{equation}
When $\eta=0$ the  contour consists of two halflines $v=\pm \sqrt{3}\,u$, $u\ge0$. For $\eta=1$ and $\eta=0$  the contours are  shown in Figure~\ref{fig:fig01}.

When $\eta<0$ we no longer have $\Im\phi\left(\sqrt\eta\right)=0$, and we need to consider both saddle points on the imaginary axis.  For the saddle point on the positive imaginary axis we write  $\beta=\sqrt{-\eta}$, we see that $\Im\phi(i{\beta})=\frac23\beta^{3}$, which gives for the path through $t=i{\beta}$
\begin{equation}\label{eq:intro03}
u=\left(v-{\beta}\right) \sqrt{\frac{v+2{\beta}}{3v}}, \quad v>0, \quad \beta=\sqrt{-\eta}.
\end{equation}
There is a similar path through $t=-i{\beta}$, and the original contour $\calC$ is split up into two parts. See  Figure~\ref{fig:fig01} for the case $\eta=-1$. Of course, we can take both paths if $f(t)$ is analytic in a domain around the contours.

When we take $f(t)=1$ the integral in \eqref{eq:intro01} becomes one of the Airy functions, in that case we have
\begin{equation}\label{eq:intro04}
F_{f=1}(\eta)=\Ai\left(\eta\right).
\end{equation}
This function satisfies the differential equation $w^{\prime\prime}(z)-zw(z)=0$, and this is the simplest second-order linear differential equation with a {\em turning point} (a point where the character of the solutions changes from oscillatory to exponential).
The Airy function $\Ai(z)$ is exponentially decaying for $z >0$, and oscillating for $z <0$.

Another classical form is
\begin{equation}\label{eq:intro05}
G(\eta)=\frac{1}{2\pi i}\int_{-\infty}^\infty e^{i\left(\frac13s^3-\eta s\right)}f(s)\,ds,
\end{equation}
which emphasizes the oscillatory nature of the integral. When $f(s)=1$ we have $G_{f=1}(\eta)=-i\Ai(-\eta)$. 

The integral in \eqref{eq:intro05} can be written (when the analytic properties of $f(s)$ allow) in the form
\begin{equation}\label{eq:intro06}
G(\eta)=\frac{1}{2\pi i}\int_{\infty e^{\frac{5}{6}\pi i}}^{\infty e^{\frac{1}{6}\pi i}} e^{i\left(\frac13s^3-\eta s\right)}f(s)\,ds,
\end{equation}
and now the oscillations at infinity are under control because the contour of integration starts and terminates in the valleys of the integrand. By substituting $s=it$, we find
\begin{equation}\label{eq:intro07}
G(\eta)=-\frac{1}{2\pi }\int_{\calC} e^{\frac13t^3+\eta t}f(it)\,dt,
\end{equation}
where $\calC$ is the contour used in \eqref{eq:intro01}.

The integrals in \eqref{eq:intro01} and \eqref{eq:intro05} arise in asymptotic analysis when we transform a contour integral like
\begin{equation}\label{eq:intro08}
H(\lambda,\alpha)=\frac{1}{2\pi i}\int_{\calL} e^{\lambda p(s,\alpha)}q(s)\,ds,
\end{equation}
where $\lambda$ is a large parameter, and when the saddle points (the zeros of $\frac{d}{ds}p(s,\alpha)$) are close together, or coalesce, for example under the influence of the parameter $\alpha$.
These integrals occur frequently in the asymptotic analysis of special functions, but also in problems from physics. 

In \S\ref{sec:transBes} we show how to transform an integral for $J_\nu(\nu z)$ into our standard form \eqref{eq:intro01}.  The function $J_\nu(\nu z)$  has a turning point at $z=1$. To see this, for the function  $J_\nu\left(\nu e^z\right)$ we have  the differential equation (see \cite[\S10.13]{Olver:2010:BFS}) 
\begin{equation}\label{eq:intro09}
w^{\prime\prime}(z)-\nu^2\left(e^{2z}-1\right)w(z)=0,
\end{equation}
with a turning point at $z=0$. Indeed, when $\nu$ is large, the function $J_\nu(\nu z)$ oscillates strongly for $z>1$ and is exponentially small for $0<z<1$.

In \S\ref{sec:transHerm} we give an  example for  the Hermite polynomial. We consider  the form $H_n\left(\sqrt{2n+1}\,x\right)$, which has Airy-type behavior when $x$ crosses the two turning points $\pm1$. Details on these examples can be found in \cite[Chapter 23]{Temme:2015:AMI}, where the transformation of  integrals of these functions into an Airy-type form is described. 
 
For Airy-type integrals of the form \eqref{eq:intro01} asymptotic expansions can be obtained. These are well defined in the neighborhood of the turning point  ($z=1$ in the case of the Bessel function $J_\nu(\nu z)$). Usually, the coefficients of the asymptotic series are difficult to evaluate near the turning point. 
There are several methods to avoid these problems, for example for the Bessel function by expanding the coefficients in powers of $(z-1)$. See \cite[\S8.5]{Gil:2007:NSF} and \cite{Temme:1997:NAU}. In \cite{Dunster:2017:CAE} a new approach is described based on computing Cauchy's integral formula for coefficients of the expansion by using the trapezoidal rule.

In this paper we show how to use a simple  numerical  method to compute  the Airy-type integrals  of the form \eqref{eq:intro01}. We use the ideas of asymptotic analysis to obtain suitable contours, the saddle point contours, without deriving asymptotic expansions.   We use the trapezoidal rule for  the integral in \eqref{eq:intro01} for  three cases: $\eta>1$, $-1\le\eta\le1$, and $\eta<-1$.  In \S\ref{sec:trapalongR} we give a few details about the trapezoidal rule for integrals on $\RR$. In a special case we compare the trapezoidal rule with Hermite quadrature. For complex values of $\eta$ the method can be used as well. Special care needs the selection of smooth contours when the phase of $\eta$ becomes $\pm\frac23\pi$, which is related with the Stokes phenomenon for the Airy-type type integral. For details we refer to \cite{Gil:2002:CCA}, where the computation of the Airy function $\Ai(z)$ is considered for complex values of $z$.
 
Apart from using  the simple trapezoidal rule, another benefit of our method is that we can use it for integrals of the form in \eqref{eq:intro08} without transforming it into an Airy-type form as in \eqref{eq:intro01}.
In \S\ref{sec:org} we use an integral representation of  the $J$-Bessel function to show how to avoid the cubic transformation.

In a recent paper \cite{Huybrechs:2018:ANM}, the approach also avoids developing an asymptotic expansion and computing the coefficients of the Airy-type asymptotic expansion and of the Airy functions. These authors use also   numerical computation of the integral in the representation given in \eqref{eq:intro01}. Their novel idea is based on Gauss quadrature on the complex contour $\calC$ with polynomials that are orthogonal on this contour. This requires the computation of the zeros of the polynomials and moments, and these topics are discussed in detail in  \cite{Huybrechs:2018:ANM}; an earlier paper on this topic is \cite{Deano:2009:CGQ}. In this way a uniform method of computation is obtained valid for general and certainly small values of the parameter $\eta$ in \eqref{eq:intro01}. In  \cite{Huybrechs:2018:ANM} an asymptotic estimate is given for error terms in their Gauss quadrature approach. We see in that paper, once the nodes and weights of the quadrature rule have been made available, very good performances of the quadrature rule for rather small degree of the underlying orthogonal polynomials.

There is a vast literature on oscillating integrals, old and new, where integrals are discussed that are more general than the Airy-type integral considered here, with asymptotic aspects and numerical evaluations. We mention research on the Pearcey integral, the swallow tail integral, with applications in physics, optics, chemical physics,  and so on. For further information, we refer to
 \cite{Borghi:2016:COS}, \cite{Connor:1990:PMU}, \cite{Connor:1992:UAO}, \cite{Deano:2018:CHO},
  \cite{Ferreira:2018:AES}, \cite{Kirk:2000:NEO}, \cite{Lopez:2017:AFE},  \cite{Milovanovic:2017:CIH}, and to references in these publications. In \cite[\S36.14]{Berry:2010:IWC} an overview of applications is given with references. For information on the Airy and Bessel functions used in this paper we refer to  
  \cite{Olver:2010:ARF} and \cite{Olver:2010:BFS}.

\section{The trapezoidal rule on \protectbold{\RR}}\label{sec:trapalongR}
The integral 
\begin{equation}\label{eq:trap01}
F(\lambda)=\intr e^{-\lambda t^2} f(t)\,dt,\quad \lambda>0,
\end{equation}
arises in asymptotic analysis as a standard form  after transformations when using Laplace's method (see, for example,  \cite[Chapter 3]{Temme:2015:AMI}). The saddle point is located at the origin and the simple dominant exponential $e^{-\lambda t^2}$ of the integrand is a convenient starting point to obtain asymptotic expansions. At the end of this section we explain how this works. After a substitution to get  $\lambda$ out of the exponential function (see \eqref{eq:trap05}), the  integral in \eqref{eq:trap01} becomes the well-known form to apply Gauss-Hermite quadrature.
We will obtain  \eqref{eq:trap01}  also after some transformations applied to our standard integral in \eqref{eq:intro01}.

We assume that $f(t)$ is analytic inside a strip of width $2a$ around $\RR$, $a>0$, and for large $t$ of order $\bigO\left(t^\alpha\right)$ for some $\alpha$. Then the trapezoidal rule has the form
\begin{equation}\label{eq:trap02}
F(\lambda)=h\sum_{k=-\infty}^\infty e^{-\lambda (kh)^2}f(hk)+R(h), \quad h >0.
\end{equation}
By using contour integration along the boundaries of the strip, using residue calculus, and assuming some decay of $\vert f(x\pm ia)\vert$ for large $\vert x\vert$, it can be shown (see  \cite{Goodwin:1949:EIF} and  \cite[\S5.4]{Gil:2007:NSF}) that the error term $R(h)$ can be estimated by
\begin{equation}\label{eq:trap03}
\vert R(h)\vert\le\frac{e^{-\pi a/h+\lambda a^2}}{\sinh(\pi a/ h)}\intr e^{-\lambda x^2} \left(\vert f(x + ia)\vert+\vert f(x - ia)\vert\right)\,dx.
\end{equation}

We combine the exponential term in front of the integral with the exponential behavior of the hyperbolic  function, and consider  $-2\pi a/h+\lambda a^2$ as a function of $a$, which is minimal for $a=\pi/(\lambda h)$. Assuming that this value of $a$ is possible (which depends on the singularities of  the function $f(t)$), and assuming that the integral in \eqref{eq:trap03} is not a large factor, we can conclude that 
\begin{equation}\label{eq:trap04}
\vert R(h)\vert=\bigO\left(e^{-\pi^2/(\lambda h^2)}\right), \quad h\to 0.
\end{equation}
This exponentially small behavior of the error was first observed by \cite{Goodwin:1949:EIF}. More details can be found in \cite[\S5.4]{Gil:2007:NSF}. An extensive recent overview on the exponentially convergent trapezoidal rule can be found in \cite{Trefethen:2014:TEC}.

In the applications from asymptotics, the parameter $\lambda$ may be large. In that case the integral in \eqref{eq:trap01} converges very fast and we need small values of $h$ such that $\lambda h^2\to 0$.  However, the number of terms needed  for numerical convergence of the series in \eqref{eq:trap02} may be acceptable. This can be seen when 
we rewrite the integral \eqref{eq:trap01}  in the form
\begin{equation}\label{eq:trap05}
F(\lambda)=\frac{1}{\sqrt{\lambda}}\intr e^{-s^2} f\left(s/\sqrt{\lambda}\right)\,ds.
\end{equation}
For this integral the estimate in \eqref{eq:trap04} has $\lambda=1$. For large values of $\lambda$ function values of  $f$ are needed for small values of its argument, as also follows from \eqref{eq:trap01}. For the algorithm there is no need to use the integral in \eqref{eq:trap05}, but it may be convenient to work always with the Gaussian function $e^{-s^2}$. 
 
 On the other hand, when in \eqref{eq:trap01} $\lambda$ is large, a few coefficients $f_k$ of the expansion $\dsp{f(t)=\sum_{k=0}^\infty f_k t^k}$ substituted in \eqref{eq:trap01} may produce the desired numerical results from the asymptotic expansion. The expansion is, under some conditions on $f(t)$,
 \begin{equation}\label{eq:trap06}
F(\lambda)\sim \sum_{k=0}^\infty f_k \intr e^{-\lambda t^2}t^k\,dt,
\end{equation}
and evaluating the integrals, observing that terms with odd index $k$ do not contribute, we obtain the classical result
 \begin{equation}\label{eq:trap07}
F(\lambda)\sim \sqrt{\frac{\pi}{\lambda}}\ \sum_{k=0}^\infty f_{2k} \left(\tfrac12\right)_k\frac{1}{\lambda^k},\quad \lambda\to\infty,
\end{equation}
where $(a)_k=\Gamma(a+k)/\Gamma(a)$ is the Pochhammer symbol.

We can efficiently use the trapezoidal rule on finite interval when the integrand is smooth with many of its successive derivatives, which are equal at the endpoints. This follows from Euler's summation formula. For an example, see the integral of the Bessel function $J_\nu(\nu z)$ on the interval $[-\pi,\pi]$ in \eqref{eq:Jbesor03}. For more details we refer to \cite[\S5.2.3]{Gil:2007:NSF}.

 \begin{table}
\caption{
Relative errors in the computation of three functions by using Gauss-Hermite  quadrature and the trapezoidal rule.
\label{tab:tab01}}
$$
\begin{array}{lccc}
{\rm case} &\ \  {\rm error\ G-H}\ \ &\ \   {\rm error \ Trap.}\ \ & \tfrac12h  \\
\hline
f(t)=1 & 2.00{\rm e-}15 & 1.00{\rm e-}15  & 1.00{\rm e-}15\\
f(t)=\cos(4t) &  2.71{\rm e-}12 &  1.09{\rm e-}06 &  1.00{\rm e-}15 \\
f(t)=1/(1+t^2) &  1.01{\rm e-}05 & 5.12{\rm e-}05 & 2.06{\rm e-}10 \\
\hline
\end{array}
$$
\end{table}

It is one of the convenient properties of the trapezoidal rule that when halving the stepsize $h$ the previous function values can be reused.  In fact,  the trapezoidal rule can be implemented as a recursive process for which there is control of the convergence rate, while Gaussian quadrature is less flexible because computing nodes and weights is more expensive and they can not be reused when the number of nodes is increased. Therefore, either a good a-priori estimation of the error is available, or one just have to just try the number of Gaussian nodes.

\begin{example}\label{ex:ex01}
In Table~\ref{tab:tab01} we compare the performance of Gauss-Hermite  with the trapezoidal rule. Gauss-Hermite should be exact for $f(t)=1$; the result gives a check that the weights and nodes are correct. 
We give three examples for the computation of the integral in \eqref{eq:trap01} with $\lambda=1$. We use Gauss-Hermite quadrature with $n=24$. For this we need 12 function evaluations, because of the symmetry of the rule and the chosen even functions. For the trapezoidal rule we compute $t_0$ such that $e^{-t_0^2}=10^{-16}$, that is, $t_0=6.0697\ldots$, and we choose $h=t_0/12=0.5058\ldots$. Hence, we use 13  nodes $k h$, $k=0,1,2,\ldots, 12$. For the results in Table~\ref{tab:tab01}  we use Maple with Digits = 16. For the numerical examples given in the paper no symbolic properties of Maple have been used.
 
We see that for the function $f(t)=\cos(4t)$ the performance of the Gauss-Hermite rule is much better. In the final column we use the trapezoidal rule with $h/2$ and 25 nodes. 
 \eoexample
\end{example}

\section{Preparing Airy-type integrals for the trapezoidal rule}\label{sec:airytype}
To evaluate the integral in \eqref{eq:intro01} by using the trapezoidal rule we consider 
  three cases: $\eta>1$, $-1\le\eta\le1$, and $\eta<-1$. This splitting up of $\RR$ is quite convenient, but we could have chosen other  intervals. The main requirement is that we do not use the saddle point contour for small values of $\vert\eta\vert$. 
  
 As explained in \S\ref{sec:intro}, we take different intervals because when $\eta\to0$ the ideal saddle point contour for numerical quadrature becomes non-smooth,  with as limiting form the two lines in Figure~\ref{fig:fig01}. Apart from this, the functions in the integrand that arise during the transformations of the variables will have their singularities approaching the origin as $\eta\to0$.

As will appear, we use transformations that will give new integrals that are special cases of those for the $K$- and $J$-Bessel functions of order $\pm\frac13$.

\subsection{The trapezoidal rule for \protectbold{\eta>1}}\label{sec:etalarge}
After the substitution $t=\sqrt{\eta}\,w$ \eqref{eq:intro01} becomes
\begin{equation}\label{eq:etalarge01}
F(\eta)=\frac{\sqrt{\eta}}{2\pi i}\int_{\calC}e^{-\frac32\xi\phi(w)}g(w)\,dw, \quad 
            \phi(w)=-\tfrac13w^3+w,
\end{equation}
where 
\begin{equation}\label{eq:etalarge02}
 \xi=\tfrac23\eta^{\frac32}, \quad g(w)=f\left(\sqrt{\eta}\,w\right).
\end{equation}

We write $w=u+iv$, which gives 
\begin{equation}\label{eq:etalarge03}
\begin{array}{@{}r@{\;}c@{\;}l@{}}
\phi(w)&=&\phi_r(u,v)+i\,\phi_i(u,v),\\[8pt]
\phi_r(u,v)&=&-\frac13u^3+uv^2+ u,\\[8pt]
\phi_i(u,v)&=&-(u^2-\frac13v^2-1)v. 
\end{array}
\end{equation}
The saddle point contour through the saddle point $w=1$ follows from solving $\phi_i(u,v)=\phi_i(1,0)=0$, that is, when we take $u=\sqrt{\frac13v^2+1}$, $v\in\RR$. We write
\begin{equation}\label{eq:etalarge04}
u=\cosh\left(\tfrac13\theta\right),\quad v=\sqrt{3}\,\sinh\left(\tfrac13\theta\right),\quad \theta\in\RR,
\end{equation}
which is a parametrization of the saddle point contour. It follows that 
\begin{equation}\label{eq:etalarge05}
\tfrac32\phi_r(u,v)=4\cosh^3\left(\tfrac13\theta\right)-3\cosh\left(\tfrac13\theta\right)=\cosh\theta,\quad \phi_i(u,v)=0.
\end{equation}
We integrate with respect to $\theta$ and obtain with $w=u+iv$ and $u,v$ as in \eqref{eq:etalarge04},
\begin{equation}\label{eq:etalarge06}
F(\eta)=\frac{e^{-\xi}}{2\pi }\sqrt{\tfrac{1}{3}\eta}\intr e^{-\xi(\cosh(\theta)-1)}g(w)\cosh\left(\tfrac13\theta\right)
\left(1-i\frac{v}{3u}\right)\,d\theta.
\end{equation}

To obtain the form in \eqref{eq:trap01}, we substitute $\tau=2\sinh\left(\frac12\theta\right)$, and we have\begin{equation}\label{eq:etalarge07}
F(\eta)=\frac{e^{-\xi}}{2\pi}\sqrt{\tfrac13\eta}\intr e^{-\frac12\xi\tau^2}h(\tau)\,d\tau,
\end{equation}
where
\begin{equation}\label{eq:etalarge08}
h(\tau)=g(w)\left(1-i\frac{v}{3u}\right)\frac{\cosh\left(\frac13\theta\right)}{\cosh\left(\frac12\theta\right)}.
\end{equation}

There are several steps needed to compute $h(\tau)$ from \eqref{eq:etalarge08}. We summarize:
\begin{enumerate}
\item
given $\tau$, compute $\theta=2\arcsinh(\frac12\tau)$;
\item
compute $w=u+iv =\cosh\left(\tfrac13\theta\right)+i\sqrt{3}\,\sinh\left(\tfrac13\theta\right)$;
\item
compute $g(w)=f\left(\sqrt\eta\,w\right)$; see \eqref{eq:etalarge02};
\item
with these values compute  $h(\tau)$.
\end{enumerate}

An important feature of the representations in \eqref{eq:etalarge06} and \eqref{eq:etalarge07} is $e^{-\xi}$ in front of the integrals. When $\eta$ is large, it is convenient to have this dominant factor explicitly outside the integral. This factor arises when we integrate the integral in \eqref{eq:etalarge01} through the saddle point  at $w=1$. On the other hand, when $\eta\to 0$, the parameter $\xi=\frac23\eta^{3/2}$ tends to 0 as well, the exponential function in the integrand loses its dominant role, and  the integrals in  \eqref{eq:etalarge01} and  \eqref{eq:etalarge07} may become divergent, due to transformations. Therefore, we consider $-1\le \eta\le1$ as a separate case. 

\begin{remark}\label{rem:rem01}
When we take in \eqref{eq:etalarge06}  $g(w)=1$, and observe that $v/u$ is odd, we obtain a $K$-Bessel function:
\begin{equation}\label{eq:etalarge09}
F(\eta)=\frac{1}{2\pi}\sqrt{\tfrac{1}{3}\eta}\intr e^{-\xi\cosh \theta}\cosh\left(\tfrac13\theta\right)\,d\theta=\frac{1}{\pi }\sqrt{\tfrac{1}{3}\eta}\,K_{\frac13}(\xi).
\end{equation}
Because of the well-known relation between this Bessel function and the Airy function (see \cite[\S9.6(i)]{Olver:2010:ARF}), we have
\begin{equation}\label{eq:etalarge10}
F(\eta)=\Ai\left(\eta\right),
\end{equation}
which also follows from \eqref{eq:intro01}  (with $f(t)=1$) and \eqref{eq:intro05}.
\eoremark
\end{remark}

\begin{remark}\label{rem:rem02}
When $\xi$ is large, we can expand the function $h(\tau)$ in \eqref{eq:etalarge08} in powers of $\tau$ and obtain an asymptotic expansion of $F(\eta)$. From the expansion $\dsp{h(\tau)=\sum_{k=0}^\infty h_k\tau^k}$ we obtain (see also \eqref{eq:trap04} and \eqref{eq:trap05})
\begin{equation}\label{eq:etalarge11}
F(\eta)\sim\frac{e^{-\xi}}{\sqrt{2\pi\xi}}\sqrt{\tfrac{1}{3}\eta}\sum_{k=0}^\infty h_{2k}2^k\left(\tfrac12\right)_k\frac{1}{\xi^k},\quad \xi\to\infty. 
\end{equation}
Similar as described after \eqref{eq:etalarge08} on computing $h(\tau)$, we need several steps to obtain the coefficients $h_k$. First,  by using \eqref{eq:etalarge04}, expand $w=u+i v =1+\sum_{k=1}^\infty w_k\theta^k$. Next expand $w$ in powers of $\tau$ by using $\theta=2\arcsinh(\frac12\tau)$. After this, we need to expand the factors of $h(\tau)$ shown in  \eqref{eq:etalarge08} in powers of $\tau$.
\eoremark
\end{remark}

\subsection{The trapezoidal rule for \protectbold{\eta<-1}}\label{sec:etaneg}
The path of integration $\calC$ in \eqref{eq:intro01} is split up into two parts $\calC=\calC_-\cup \calC_+$,  with $\calC_-$ from $\infty e^{-\frac13\pi i}$ to $-\infty$, and $\calC_+$ from $-\infty$ to $\infty e^{\frac13\pi i}$.  See Figure~\ref{fig:fig01}, where we see the two parts for the case $\eta=-1$ running from  $-\infty$ to $\infty e^{\pm\frac13\pi i}$, with the indicating directions of integration. As long as $\vert\eta\vert$ remains bounded away from zero, the contours $\calC_\pm$ are smooth, just as the original contour $\calC$ in the previous  case for $\eta>1$.

We denote the contributions from the paths $C_\pm$ by $F^{\pm}(\eta)$, and we have $F(\eta)=F^+(\eta)+F^-(\eta)$. The saddle points are at $\pm i\beta$, with $\beta=\sqrt{-\eta}$, and the steepest descent paths run through these points. We assume that $f(t)$ is analytic around these paths and that for large values of $t$ we have $f(t)=\bigO\left(t^\alpha\right)$ for some $\alpha$.

When we have computed $F^+(\eta)$ and $f(t)$ is a real function (real for real values of $t$), then we find $F(\eta)$ by taking twice the real part of $F^+(\eta)$. That is, $F(\eta)=2\Re F^+(\eta)$.

We consider the contribution  from the path in the upper half plane
\begin{equation}\label{eq:etaneg01}
F^+(\eta)=\frac{1}{2\pi i}\int_{\calC_+}e^{\frac13t^3+\beta^2 t}f(t)\,dt,\quad \beta=\sqrt{-\eta},\quad \beta>1,
\end{equation}
through the saddle point $t=+i\beta$. We write $t=\beta w$, which gives
\begin{equation}\label{eq:etaneg02}
F^+(\eta)=\frac{\beta}{2\pi i}e^{i\xi}\int_{\calC_+}e^{-\frac32\xi \phi(w)}g(w)\,dw,
\end{equation}
where
\begin{equation}\label{eq:etaneg03}
g(w)=f(\beta w), \quad \phi(w)=-\tfrac13w^3-w+\tfrac23i,
\quad \xi=\tfrac23\beta^{3}.
\end{equation}

At the saddle point $w=i$ we have $\phi(i)=0$. With $w=u+iv$, the saddle point contour follows from the equation 
$\Im\phi(w)=-u^2v+\frac13v^3-v+\frac23=0$. A  solution is  $u=(v-1)\sqrt{(v+2)/(3v)}$, $v>0$. In the Introduction we have obtained this form (see \eqref{eq:intro03}), but here we used the substitution $t=\beta w$.

This time we use a different parametrization of the saddle point contour, because we like to obtain an integral from which a  representation of the $J$-Bessel function follows, in a similar manner as we have obtained \eqref{eq:etalarge06}, which becomes a $K$-Bessel function when $g(w)=1$. See also Remark~\ref{rem:rem01}.

\begin{figure}
\begin{center}
\epsfxsize=10cm \epsfbox{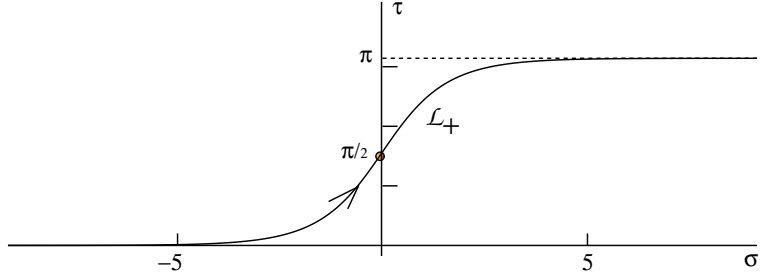}
\caption{
\label{fig:fig02} Saddle point contour $\calL_+$ defined by \eqref{eq:etaneg06} in the $\theta$-plane, $\theta=\sigma+i\tau$, $\tau\in(0,\pi)$.  }
\end{center}
\end{figure}

We take in \eqref{eq:etaneg01} $w=2\sinh\left(\frac13\theta\right)$, write $\theta=\sigma+i\tau$, and obtain $\phi(w)=\phi_r(\sigma,\tau)+i\phi_i(\sigma,\tau)$, where
\begin{equation}\label{eq:etaneg04}
\begin{array}{@{}r@{\;}c@{\;}l@{}}
\phi(w)&=&-\frac23\sinh\theta+\frac23i,\\[8pt]
\phi_r(\sigma,\tau)&=&-\frac23\sinh\sigma\cos\tau,\\[8pt]
\phi_i(\sigma,\tau)&=&\frac23\cosh\sigma\sin\tau-\frac23,
\end{array}
\end{equation}
and
\begin{equation}\label{eq:etaneg05}
F^+(\eta)=\frac{\beta}{3\pi i}e^{i\xi}\int_{\calL_+}e^{\xi \left(\sinh\theta-i\right)}g(w)\cosh\left(\tfrac13\theta\right)\,d\theta.
\end{equation}
The path $\calL_+$ is shown in Figure~\ref{fig:fig02} and follows from $\phi_i(u,v)=0$. It is given by
\begin{equation}\label{eq:etaneg06}
\cosh \sigma \sin \tau=1,\quad 0 < \tau < \pi.
\end{equation}
When we integrate with respect to $\sigma$, we need ${d\tau}/{d\sigma}$. From \eqref{eq:etaneg06} we obtain for $\sigma>0$, where $\tau\in(\frac12\pi,\pi)$ (see Figure~\ref{fig:fig02}), and, hence, $\cos\tau <0$,
\begin{equation}\label{eq:etaneg07}
\frac{d\tau}{d\sigma}=-\frac{\sinh\sigma}{\cos\tau\,\cosh^2\sigma}=\frac{\sinh\sigma}{\sqrt{1-\sin^2\tau}\,\cosh^2\sigma}=\frac{1}{\cosh\sigma}=\sin\tau.
\end{equation}
When $\sigma<0$ we have $\cos\tau>0$ and $\sinh\sigma<0$, and we obtain the same result for ${d\tau}/{d\sigma}$. 

It follows that
\begin{equation}\label{eq:etaneg08}
F^+(\eta)=\frac{\beta}{3\pi i}e^{i\xi}\int_{-\infty}^\infty e^{\xi\sinh\sigma\cos\tau}g(w)\cosh\left(\tfrac13\theta\right)(1+i\sin\tau) \,d\sigma.
\end{equation}

For the argument of the exponential function we have  $\sinh\sigma\cos\tau\le0$ for all $\sigma\in\RR$. To see this, 
we have using \eqref{eq:etaneg06}
\begin{equation}\label{eq:etaneg09}
\cos^2\tau=\frac{\sinh^2\sigma}{\cosh^2\sigma}\quad \Longrightarrow \quad \cos \tau=-\frac{\sinh\sigma}{\cosh\sigma},
\end{equation}
where the minus sign is chosen because  $\sign(\cos\tau)=-\sign(\sigma)$. In this way, we obtain 
\begin{equation}\label{eq:etaneg10}
\sinh\sigma\cos\tau=-\frac{\sinh^2\sigma}{\cosh\sigma}=-\tanh\sigma\sinh\sigma,\quad \sigma\in\RR
\end{equation}
and
\begin{equation}\label{eq:etaneg11}
F^+(\eta)=\frac{\beta}{2\pi i}e^{i\xi}\int_{-\infty}^\infty e^{-\xi\tanh\sigma\sinh\sigma}g(w)\cosh\left(\tfrac13\theta\right)(1+i\sin\tau) \,d\sigma.
\end{equation}

\begin{remark}\label{rem:rem03}
When  we take  $g(w)=1$, the  integral in \eqref{eq:etaneg05} reduces to one that defines a sum of $J$-and $Y$-Bessel functions of order $\pm\frac13$, similarly as  in Remark~\ref{rem:rem01}. To verify this, we have for the $J$-Bessel function the representation
\begin{equation}\label{eq:etaneg12}
J_\nu(z)=\frac{1}{2\pi i}\int_{\calL}e^{z\sinh s-\nu s}\,ds.
\end{equation}
The contour $\calL$  in \eqref{eq:etaneg12} starts at $+\infty-\pi i$ and terminates at $+\infty+\pi i$.  It may visit $-\infty$ on the negative real axis, and then it can be composed by using the contours for the Hankel functions. In this way,
we obtain
\begin{equation}\label{eq:etaneg13}
\begin{array}{@{}r@{\;}c@{\;}l@{}}
H_\nu^{(1)}(z)&=&\dsp{\frac{1}{\pi i}\int_{-\infty}^{\infty +\pi i}e^{z\sinh s-\nu s}\,ds,}\\[8pt]
H_\nu^{(2)}(z)&=&\dsp{\frac{-1}{\pi i}\int_{-\infty}^{\infty -\pi i}e^{z\sinh s-\nu s}\,ds.}
\end{array}
\end{equation}
Together these functions give
\begin{equation}\label{eq:etaneg14}
J_\nu( z)=\tfrac{1}{2}\left(H_\nu^{(1)}(z)+H_\nu^{(2)}(z)\right).
\end{equation}

It follows that, when $g(w)=1$, \eqref{eq:etaneg05} becomes
\begin{equation}\label{eq:etaneg15}
\begin{array}{@{}r@{\;}c@{\;}l@{}}
F^+(\eta)&=&\frac16 \beta\left(H_{-\frac13}^{(1)}(\xi)+H_{\frac13}^{(1)}(\xi)\right)\\[8pt]
&=&\frac16\beta \left(J_{-\frac13}(\xi)+iY_{-\frac13}(\xi)+J_{\frac13}(\xi)+iY_{\frac13}(\xi)\right).
\end{array}
\end{equation}
As remarked at the beginning of this section, the complete  $F(\eta)$ follows from taking twice the real part of $F^+(\eta)$, and we obtain
\begin{equation}\label{eq:etaneg16}
F(\eta)=\tfrac13\beta\left(J_{-\frac13}(\xi)+J_{\frac13}(\xi)\right),\quad  \beta=\sqrt{-\eta}, \quad \xi=\tfrac23\beta^3.
\end{equation}
This corresponds with \cite[Eqn.~9.6.6]{Olver:2010:ARF}
\begin{equation}\label{eq:etaneg17}
\Ai(-z)=\tfrac13\sqrt{z}\left(J_{-\frac13}(\zeta)+J_{\frac13}(\zeta)\right),\quad  \zeta=\tfrac23z^{\frac32}.
\end{equation}
\eoremark
\end{remark}

\subsection{The trapezoidal rule for \protectbold{-1\le\eta\le1}}\label{sec:etasmall}

For large positive values of $\eta$ the function $F(\eta)$ defined in \eqref{eq:intro01} will be very small. The front factor $e^{-\xi}$ in \eqref{eq:etalarge01} and in later formulas describes this behavior very well. This factor arises when we choose for the contour $\calC$ in \eqref{eq:intro01} or in \eqref{eq:etalarge01} a contour through the saddle points $t=\sqrt{\eta}$ or $w=1$.

However, because of the transformation $t=\sqrt{\eta}\,w$, the integrals in \eqref{eq:etalarge06}  and  \eqref{eq:etalarge07} become useless when $\eta\to0$ and, hence, $\xi\to0$. To handle this for the present
values of $\eta$, we do not use the $t\to w$ transformation and we do not use a path through the saddle point    $t=\sqrt{\eta}$. 
As a consequence,  we miss the  factor $e^{-\xi}$ as in \eqref{eq:etalarge04}, which is dominant when $\xi$ is large, but this factor is not relevant in the present case. 

We take for the path of  the integral in \eqref{eq:intro01}  $t=u+iv$, with $u=1+\sqrt{1+\frac13v^2}$. 
In this way the path is independent of $\eta$, it cuts the real axis at the fixed point $t=2$, and it runs into the valleys of the function $e^{\frac13t^3}$. We could have taken $u$ slightly different, but the present choice works well and is convenient.

We parametrize the path by writing
\begin{equation}\label{eq:etasmall01}
u=1+\cosh\theta,\quad v=\sqrt{3}\,\sinh\theta,\quad \theta\in\RR.
\end{equation}
This gives
\begin{equation}\label{eq:etasmall02}
\begin{array}{@{}r@{\;}c@{\;}l@{}}
\phi(t)&=&\dsp{\tfrac13t^3-\eta t=\tfrac83-2\eta-p(\theta)+ir(\theta),}\\[8pt]
p(\theta)&=&\frac13(\cosh\theta-1)\left(8\cosh^2\theta+14\cosh\theta+2+3\eta\right),\\[8pt]
r(\theta)&=&\sqrt{3}\,\sinh\theta\left(2\cosh\theta + 2-\eta\right).
\end{array}
\end{equation}

When we  integrate  with respect to $\theta$, using $\dsp{\frac{dt}{d\theta}=\sinh\theta+i\sqrt3\cosh\theta}$,
we obtain 
\begin{equation}\label{eq:etasmall03}
F(\eta)=\frac{\sqrt{3}\,e^{\frac83-2\eta}}{2\pi i }\intr e^{-p(\theta)}q(\theta)\,d\theta,\quad q(\theta)=\frac{1}{\sqrt{3}}f(t)e^{ir(\theta)}\frac{dt}{d\theta}.
\end{equation}

Because we have not chosen a saddle point contour, on our present path $\Im\phi(t)=r(\theta)$ is not constant. Apart from the choice of $f(t)$ in  \eqref{eq:intro01}, this causes oscillations due to the choice of our contour. This becomes visible in the function $q(\theta)$. Again, apart from the influence of $f(t)$, we observe that $r(\theta)$ is of lower growth than the function $p(\theta)$ for large $\theta$. The oscillations due to $e^{ir(\theta)}$ have some influence on the convergence when applying the trapezoidal rule, but the dominant factor $e^{-p(\theta)}$ will damp these oscillations quite well.

\begin{figure}
\begin{center}
\epsfxsize= 6cm \epsfbox{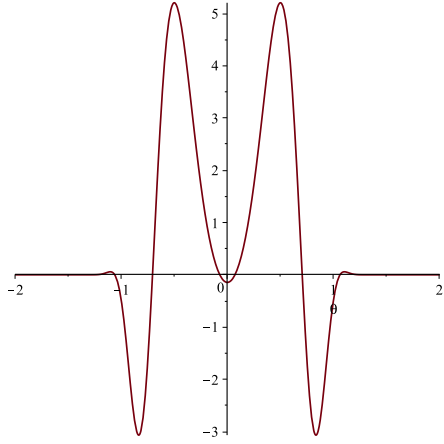}
\caption{
\label{fig:fig03} Graph of $ \Im \left(e^{-p(\theta)}q(\theta)\right)$ when $\eta=1$ and $f(t)=\cos(4t)$ in \eqref{eq:etasmall03} for $\theta\in[-2,2]$. }
\end{center}
\end{figure}

\begin{example}\label{ex:ex02}
In Figure~\ref{fig:fig03} we show a graph of the integrand part $ \Im \left(e^{-p(\theta)}q(\theta)\right)$ with $\eta=1$, $f(t)=\cos(4t)$, $\theta\in[-2,2]$. We observe some oscillations of the integrand and fast damping thereof. 
This choice of $f(t)$ causes extra oscillations and is exponentially large for large $\Im t$. When $\vert \theta\vert>1.65$ the integrand of the integral in \eqref{eq:etasmall03}
  is smaller than $1.50e-16$. We take $h=0.1$ and $k_{max}=\lfloor1.65/h\rfloor=16$, the number of terms in the trapezoidal with positive index $k$. We compare the result with the value of $\Re\left(\Ai(\eta+4i\right)$ and find a relative error $1.87e-8$. When we halve the stepsize, taking $h=0.05$ and  $k_{max}=\lfloor1.65/h\rfloor=33$, we find a relative error $3.00e-15$. These results are similar for other values of $\eta\in[-1,1]$.
\eoexample
\end{example}

\subsection{Other type of contours}\label{sec:other}
In \cite{Huybrechs:2018:ANM} the goal of the paper is the construction and analysis of a uniformly applicable quadrature rule, uniform in the parameter $c$ near $c = 0$, for the canonical integral
\begin{equation}\label{eq:other01}
I(c)=\int_{-1}^1 e^{i\omega\left(\frac13t^3-ct\right)} f(t)\,dt\approx \sum_{k=1}^n w_k f(x_k),
\end{equation}
where $f(t)$ is an analytic function at least in an open neighborhood of $[-1,1]$ in the complex plane. The interval of integration can be modified into three parts: from $t=-1$ into the valley at $\infty e^{\frac{5}{6}\pi i}$ giving an integral $I_{-1}(c)$, then from that valley into the valley at $e^{\frac{1}{6}\pi i}$, giving an integral $I_{-1,1}(c)$, and finally back to $t=1$, giving an integral $I_{1}(c)$. The integral $I_{-1,1}(c)$ is like the integral in \eqref{eq:intro06}, and this one is in \cite{Huybrechs:2018:ANM} computed with Gauss quadrature on a complex contour. Of course, $f(t)$ should be analytic in the extended domains of the complex plane.

For example, in  \cite[\S8]{Huybrechs:2018:ANM} we see the integral with the Chebyshev polynomial  
\begin{equation}\label{eq:other02}
\begin{array}{@{}r@{\;}c@{\;}l@{}}
I(k,\omega)&=&\dsp{\int_{-1}^1 T_k(x) e^{i\omega x}\,dx=\int_0^\pi \sin\theta\, \cos(k\theta)\, e^{i\omega \cos\theta}\,d\theta}\\[8pt]
&=&\dsp{\tfrac12\int_0^\pi \sin\theta\,  e^{i(k\theta+\omega \cos\theta)}\,d\theta+\tfrac12\int_0^\pi \sin\theta\, e^{i(-k\theta+\omega \cos\theta)}\,d\theta}\\[8pt]
&=&I_1(k,\omega)+I_2(k,\omega).
\end{array}
\end{equation}
The parameters $k$ and $\omega$ are large and when $k\sim\omega$ the integral $I_1(k,\omega)$ has two nearby saddle points
(we use $\mu=k/\omega$)
\begin{equation}\label{eq:other03}
\begin{array}{lll}
\theta_1=\arcsin\mu,\quad &\theta_2=\pi-\arcsin\mu,\quad &0 < \mu\le1,\\[8pt]
\theta_1=\frac12\pi+i\,\arccosh\,\mu,\quad &\theta_2=\frac12\pi-i\,\arccosh\,\mu,\quad  &\mu\ge1.\\[8pt]
\end{array}
\end{equation}

\begin{figure}
\begin{center}
\epsfxsize= 11cm \epsfbox{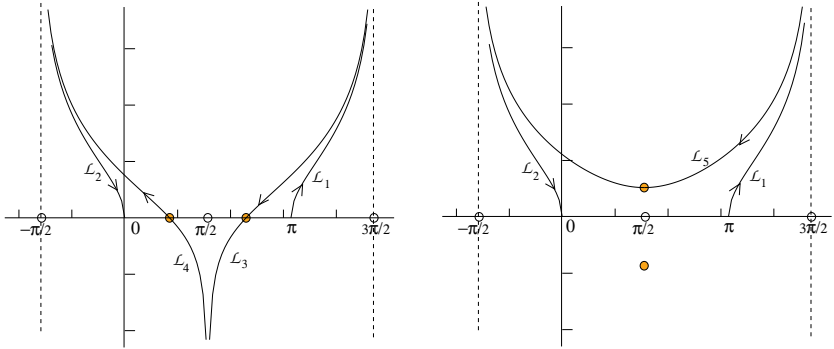}
\caption{
\label{fig:fig08} Paths on which the oscillator $e^{i(k\theta+\omega \cos\theta)}$ of the integral $I_1(k,\omega)$ given in \eqref{eq:other02}  has  constant imaginary parts. In the  figure on the left we have $\mu=k/\omega <1$, on the right $\mu=k/\omega >1$.}
\end{center}
\end{figure}

In Figure~\ref{fig:fig08} we show the paths on which the oscillator  $e^{i(k\theta+\omega \cos\theta)}$ of the integral $I_1(k,\omega)$ has  constant imaginary parts, and the integral over $[0,\pi]$ follows from the roundtrip over all contours
\begin{equation}\label{eq:other04}
\begin{array}{lll}
&[0,\pi]\,\cup\, \calL_1 \,\cup\, \calL_3\,\cup\, \calL_4\,\cup\,\calL_2,\quad &0<  \mu\le  1,\\[8pt]
&[0,\pi]\,\cup\, \calL_1 \,\cup\, \calL_5\,\cup\, \calL_2,\quad  &\mu\ge 1.\\[8pt]
\end{array}
\end{equation}
The complete integral over these 4 or 5 paths in the shown directions  is 0. The contours $\calL_3, \calL_4, \calL_5$ are saddle point contours on which we can use a transformation with a cubic polynomial to give the form of an Airy-type integral as in \eqref{eq:intro06}. However, we can also use the trapezoidal rule without this transformation, as explained for Bessel function contours in \S\ref{sec:org}.

In \cite{Gil:2001:ONI} and \cite{Gil:2002:AGH} we have considered the analysis and computation of the so-called inhomogeneous Airy functions, or Scorer functions, defined by
\begin{equation}\label{eq:other05}
\Gi(z)=\frac1{\pi}\int_0^\infty \sin\left(\tfrac13t^3+zt\right)\,dt,\quad
\Hi(z)=\frac1{\pi}\int_0^\infty e^{-\frac13t^33+zt}\,dt.
\end{equation}
The first integral is only defined for $z\in\RR$, the other one for $z\in\CC$, but we have many interrelations for these functions, such as
\begin{equation}\label{eq:other06}
\Gi(z)=-\tfrac12\left(e^{\frac23\pi i}\Hi\left(ze^{\frac23\pi i}\right)+e^{-\frac23\pi i}\Hi\left(ze^{-\frac23\pi i}\right)\right).
\end{equation}
We have constructed non-oscillating contours for the Scorer functions, and we have used the trapezoidal rule for computations.  After extra transformations, we can also use Gauss-Laguerre quadrature. 

The  integrals $I_{\pm1}(c)$ that follow from \eqref{eq:other01} can be computed in the same way, as is also proposed in  \cite{Huybrechs:2018:ANM}. We can do the integral $I_{-1,1}(c)$t with our approach, as follows from a simple transformation from \eqref{eq:intro06} to \eqref{eq:intro07}.

\section{Transforming integrals into Airy-type integrals}\label{sec:transAirytype}
We give two examples of the transformation of integrals into the standard form in \eqref{eq:intro01}, one for the $J$-Bessel function and another one for the Hermite polynomial.

\subsection{The Bessel function}\label{sec:transBes}
We use the integral representation already given in \eqref{eq:etaneg12}:
\begin{equation}\label{eq:transBes01}
J_\nu(z)=\frac{1}{2\pi i}\int_{\calL}e^{z\sinh s-\nu s}\,ds,
\end{equation}
where  we  assume that $\nu$ is positive. The contour $\calL$ starts at $+\infty-\pi i$ and terminates at $+\infty+\pi i$.  We will use this representation by replacing $z$ by $\nu z$, which gives
\begin{equation}\label{eq:transBes02}
J_\nu(\nu z)=\frac{1}{2\pi i}\int_{\calL}e^{\nu(z\sinh s- s)}\,ds,
\end{equation}
The two saddle points of this integral 
$s_\pm=\pm \arccosh (1/z)$ coalesce when $z=1$. For $0<z\le 1$ the saddle points are real, and the function $J_\nu(\nu z)$ is monotonic (as a function of $z$). 

For this integral the transformation into an Airy-type integral is originally introduced in \cite{Chester:1957:ESD} and reads
\begin{equation}\label{eq:transBes03}
z\sinh s-s=\tfrac13r^3-\zeta r+A, 
\end{equation}
where $A$ and $\zeta$ follow from substituting the corresponding saddle points $(s_+,\sqrt\zeta)$ and $(s_-,-\sqrt\zeta)$  in the $s$ and $r$ plane. This gives $A=0$ and for $\zeta$ we find
\begin{equation}\label{eq:transBes04}
\begin{array}{@{}r@{\;}c@{\;}l@{}}
\frac23\zeta^{\frac32}&=&\arccosh(1/z)-\sqrt{1-z^2},\quad 0<z\le 1,\\[8pt]
\frac23(-\zeta)^{\frac32}&=&\sqrt{z^2-1}-\arccos(1/z),\quad \quad z\ge1.
\end{array}
\end{equation}
The relation between $z$ and $\zeta$ is analytic at $z=1$, where $\zeta=0$, and we have the local expansion
\begin{equation}\label{eq:transBes05}
z(\zeta)=1-\lambda+\tfrac3{10}\lambda^2+\tfrac1{350}\lambda^3+\ldots,\quad \lambda=2^{-\frac13}\zeta.
\end{equation}
The transformation in \eqref{eq:transBes03} gives for \eqref{eq:transBes02} the standard form
\begin{equation}\label{eq:transBes06}
J_\nu(\nu z)=\frac{1}{2\pi i}\int_{\calC}e^{\nu\left(\tfrac13r^3-\zeta r\right)}h(r)\,dr,
\end{equation}
where
\begin{equation}\label{eq:transBes07}
h(r)=\frac{ds}{dr}=\frac{r^2-\zeta}{z\cosh s-1},
\end{equation}
and the contour $\calC$ runs from $\infty e^{-\frac13\pi i}$ to $\infty e^{+\frac13\pi i}$. This is the convenient contour when $\zeta\ge0$ (or $0<z\le1$), for negative values of $\zeta$ the contour can be split up as shown in Figure~\ref{fig:fig01} for $\eta=\pm1$.

We prepare the integral representation in  \eqref{eq:transBes06} for the trapezoidal rule, and substitute $r=\nu^{-\frac13}t$ . We obtain an integral as in  \eqref{eq:intro01}:
\begin{equation}\label{eq:transBes08}
J_\nu(\nu z)=\frac{1}{2\nu^{\frac13}\pi i}\int_{\calC}e^{\frac13t^3-\eta t}f(t)\,dt,\quad f(t)=h\left(\nu^{-\frac13}t\right),\quad \eta =\nu^{\frac23}\zeta.
\end{equation}
As observed earlier, when  $h(r)$ is replaced by a constant the integral becomes an Airy function. 

\begin{figure}
\begin{center}
\epsfxsize=10cm \epsfbox{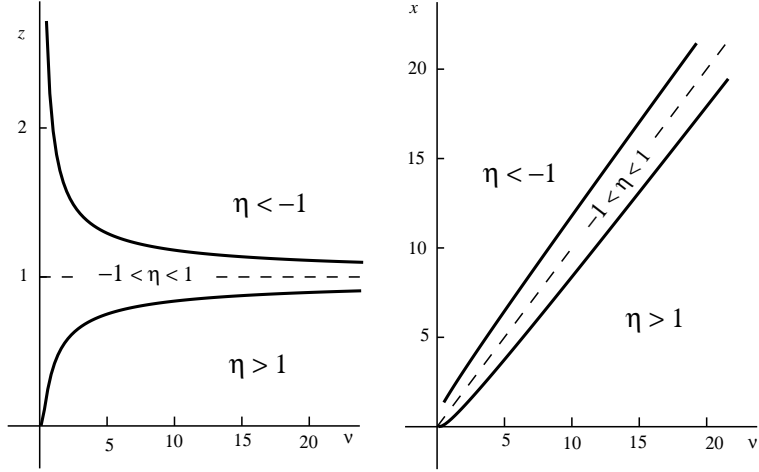}
\caption{
\label{fig:fig04} The domains in the $(\nu,z)$-plane (left, for $J_\nu(\nu z)$) and in the $(\nu,x)$-plane (right, for $J_\nu(x)$), where the three domains  $\eta < -1$, $-1 \le \eta \le 1$ and $\eta > 1$ are located. }
\end{center}
\end{figure}

Now we can consider the three cases for $\eta$ as in \S\S\ref{sec:etalarge}--\ref{sec:etasmall}.
In Figure \ref{fig:fig04} we show the domains in the $(\nu,z)$-plane (left, for $J_\nu(\nu z)$) and in the $(\nu,x)$-plane (right, for $J_\nu(x)$), where the three domains $\eta < -1$, $-1 \le \eta \le 1$ and $\eta > 1$ are located. 
We observe that for large values of $\nu$ a small interval  area for $\eta\in[-1,1]$ for $\eta$ arises.

\begin{remark}\label{rem:rem04}
Usually we start the computation of the Bessel function $J_\nu(\nu z)$ with $z$ and $\nu$ given. Then  $\zeta$ follows from  \eqref{eq:transBes04}, and \eqref{eq:transBes08} gives $\eta=\nu^{\frac23}\zeta$. So, only knowing $\nu$ and  $z$ is not enough to see in a glance which $\eta$-value arises, and which method for the trapezoidal rule follows.  Because of the simple form of the exponential function in  \eqref{eq:transBes08}, with just one parameter $\eta$, this integral representation gives a simple criterion to choose one of the three methods for the trapezoidal rule. 
\eoremark
\end{remark}

\subsection{The Hermite polynomial}\label{sec:transHerm}
The Hermite polynomials follow from the generating function
\begin{equation}\label{eq:transHerm09}
e^{2xz-z^2}=\sum_{n=0}^\infty\frac{H_n(x)}{n!}z^n,\quad x,z\in\CC,                
\end{equation}
which gives the Cauchy-type integral
\begin{equation}\label{eq:transHerm10}
H_n(x)=\frac{n!}{2\pi i}\, \int_\calC e^{2xz-z^2} \,\frac{dz}{z^{n+1}},
\end{equation}
where $\calC$ is a circle around the origin and the integration is in the
positive direction. By substituting $z=\nu s$, $x=\nu\xi$, $\nu=\sqrt{2n+1}$, we obtain
\begin{equation}\label{eq:transHerm11}
H_n(x)=\frac{n!}{\nu^n}\frac{1}{2\pi i}\, \int_\calC   e^{\nu^2\phi(s)}\,\frac{ds}{\sqrt{s}},     
\end{equation}
where 
\begin{equation}\label{eq:transHerm12}
\phi(s)=2\xi s-\tfrac12\ln s-s^2,
\end{equation}
and $\calC$ is a path that runs from $-\infty$ (with
$\phase\,s=-\pi$), encircles the origin in positive direction, and returns to
$-\infty$, now with $\phase\,s=+\pi$. If we wish we can extend the contour to $+\infty$, which we will do in the oscillatory case. 

\begin{figure}[ht]
\begin{center}
\epsfxsize=8cm \epsfbox{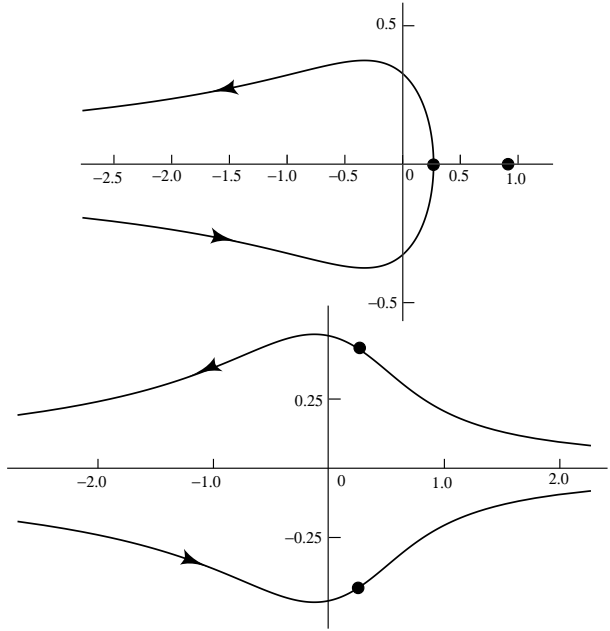}
\caption{
\label{fig:fig05} Saddle point contours for $\xi=1.2$ (upper figure, through the saddle point $s_-$) and $\xi=0.5$ (lower figure, through both saddle points). }
\end{center}
\end{figure}

We can assume that $x\ge0$ ($\xi\ge0$) because of symmetry $H_n(-x)=(-1)^nH_n(x)$, which follows from \eqref{eq:transHerm09} by changing $x\to-x$ and $z\to-z$.

The saddle points of the integral are
defined by the equation $\phi^\prime(s)=2\xi-1/(2s)-2s=0$ and are given by
\begin{equation}\label{eq:transHerm13}
s_{\pm}=\tfrac12\left(\xi\pm\sqrt{{\xi^2-1}}\right).  
\end{equation}

When $\xi=\pm1$ the saddle points coalesce at $ \frac12$, and when
$\xi\sim1$ uniform Airy-type expansions can be derived. When
$-1<\xi<1$ the saddle points are complex and are located on the circle around the origin with radius $\frac12$. For
these values of $\xi$, that is, if $-\sqrt{{2n+1}}<x<\sqrt{{2n+1}}$, zeros occur. When
$\xi>1$ or $\xi<-1$ the saddle points are real, and the Hermite polynomials are
non-oscillating.

In Figure~\ref{fig:fig05} we show the saddle point contours of the integral in \eqref{eq:transHerm11} for $\xi=1.2$ (upper figure) and $\xi=0.5$ (lower figure). The saddle points are indicated by black dots.

We can obtain an Airy-type integral  using the transformation
\begin{equation}\label{eq:transHerm14}
\phi(s)=\tfrac13r^3-\zeta r+ A,
\end{equation}
where $\phi(s)$ is defined in \eqref{eq:transHerm12}. For details we refer to  \cite[\S23.4]{Temme:2015:AMI}.

For $\zeta$ and $A$ we have the relations
\begin{equation}\label{eq:transHerm15}
\begin{array}{@{}r@{\;}c@{\;}l@{}l@{\,}}
\frac43\zeta^{\frac32}&=&\xi\sqrt{\xi^2-1}-\arccosh\,\xi, \quad &\xi\ge1,\\[8pt]
\frac43(-\zeta)^{\frac32}&=&\arccos \xi-\xi\sqrt{1-\xi^2},\quad &0\le\xi\le1,\\[8pt]
A&=&\frac12\xi^2+\frac14+\frac12\ln2.
\end{array}
\end{equation}  
These follow from substituting the corresponding saddle points in the $s$ and $r$ variables.

The quantity $\zeta$ is analytic in a neighborhood of $\xi=1$. Because of symmetry we only need $x\ge0$, that is $\xi\ge0$, 
and we have $\zeta \ge (\frac34\arccos(0))^{2/3}= -1.11546\ldots$. Also, there is simple  differential equation
\begin{equation}\label{eq:transHerm16}
\zeta\left(\frac{d\zeta}{d\xi}\right)^2=\xi^2-1,
\end{equation}
and for small values of $\vert\xi-1\vert$ there is an expansion 
\begin{equation}\label{eq:transHerm17}
\zeta=2^{\frac13}(\xi-1)\left(1+\tfrac1{10}(\xi-1)-\tfrac2{175}(\xi-1)^2+\bigO\left((\xi-1)^3\right)\right).
\end{equation}

The cubic transformation in \eqref{eq:transHerm14} gives the representation
\begin{equation}\label{eq:transHerm18}
H_n(x)=\frac{n!\,e^{\nu^2A}}{\nu^n}\frac{1}{2\pi i}\, \int_\calL   e^{\nu^2\left(\frac13r^3-\zeta r\right)}g(r)\,dr,
\end{equation}
where 
\begin{equation}\label{eq:transHerm19}
g(r)=-\frac{1}{\sqrt{s}}\frac{ds}{dr}=-\frac{1}{\sqrt{s}}\frac{\psi^\prime(r)}{\phi^\prime(s)}=2\sqrt{s}\frac{r^2-\zeta}{4s^2-4\xi s+1}.
\end{equation}
The transformation $r=t\nu^{-\frac23}$ gives
\begin{equation}\label{eq:transHerm20}
H_n(x)=\frac{n!\,e^{\nu^2A}}{\nu^{n+\frac23}}\frac{1}{2\pi i}\, \int_\calL  e^{\frac13t^3-\eta t}f(t)\,dt,\quad
\eta=\nu^{\frac43} \zeta, \quad f(t)=g\left(t\nu^{-\frac23}\right).
\end{equation}

The contour $\calL$ is a contour as  described for $\calC$ in \eqref{eq:intro01},  
see Figure~\ref{fig:fig01}, where we have shown the contour for a few values of $\eta$. For example, when $\eta>0$ in \eqref{eq:transHerm20}, the contour has the form as the one for $\eta=1$ in Figure~\ref{fig:fig01}. It is, up to scaling, the image of  the upper contour in Figure~\ref{fig:fig05} with $\xi >1$.

For $x\ge0$ the parameter $\eta$ is real. We have $\eta\ge \nu^{\frac43} 1.11546\ldots$, see above \eqref{eq:transHerm16}.

The integral representation in \eqref{eq:transHerm20}  can be  used for numerical evaluation by using the trapezoidal rule. For obtaining the Airy-type asymptotic expansion we refer to  \cite[\S23.4]{Temme:2015:AMI}.

\section{Evaluating the $J$-Bessel function by using a standard integral representation}\label{sec:org}
We have explained in \S\ref{sec:airytype} that the numerical evaluation of the Airy-type integral in \eqref{eq:intro01} is quite straightforward, but we have not considered the evaluation of the function $f(t)$. In the examples of the Bessel function and the Hermite polynomial the function arises when we use the cubic transformation; see \eqref{eq:transBes07} and \eqref{eq:transBes08} for the Bessel-case. The function $h(r)$ in \eqref{eq:transBes07} looks quite simple, but the numerical evaluation is not so easy near the saddle point $r=\sqrt{\zeta}$.
We have, by using  l'H\^{o}pital's rule,
\begin{equation}\label{eq:org01}
h\left(\sqrt{\zeta}\right)=\left(\frac{4\zeta}{1-z^2}\right)^{\frac14}.
\end{equation}

When applying the trapezoidal rule, say for $\eta \ge1$, the main contributions to the integral in \eqref{eq:transBes08} come from a small neighborhood of $t=\sqrt{\eta}$, especially for large values of $\nu$.
So, we need series expansions, for example of the form 
\begin{equation}\label{eq:org02}
f(t)=\sum_{k=0}^\infty c_k\left(t-\sqrt{\eta}\right)^k.
\end{equation}
The coefficients $c_k$ can be found in analytic form, but for small values of $\zeta$ there is another numerical problem.  For the Bessel functions, the $c_k$ are finite sums of negative powers of $\zeta$ and $(z-1)$; see  \eqref{eq:transBes04} and\eqref{eq:transBes05} for the relation between $\zeta$ and $z$. The representations of the coefficients of the Airy-type expansion of the Bessel functions are shown in in equations (12.10.10)--(12.10.13) of \cite[\S10.20(i)]{Olver:2010:BFS}. The $c_k$ have similar forms. The limit of $c_k$ as $\zeta\to0$ is well defined, because the $c_k$ are analytic in a  neighborhood of $\zeta=0$. Numerical cancellation of digits will happen because of the removable singularities at $\zeta=0$. A numerical issue already occurs in the evaluation of $h\left(\sqrt{\zeta}\right)$ in \eqref{eq:org01}: for small values of $\zeta$ we need an expansion as given in \eqref{eq:transBes05}. It is a notorious drawback of the method of uniform asymptotics that the coefficients are difficult to evaluate  when $\zeta\sim0$ (or $z\sim1$), that is, when the saddle points coalesce. 
 
For more general problems, outside the area of Hermite polynomials, Bessel functions, and other special functions, the functions arising in the integrand  from the transformation of a function into a cubic polynomial, similar difficulties will arise. From an analytical point of view, the function $f(t)$ as in \eqref{eq:transBes08} can always be expanded around the saddle point, as in  \eqref{eq:org02}, with $\eta=\nu^{\frac23}\zeta$, see \eqref{eq:transBes08}. The definition of  $\zeta$ looks always as  in \eqref{eq:transBes04}, with different right-hand sides.
Expansions like \eqref{eq:transBes05} (and inverted ones) should be derived, and next the coefficients $c_k$ of the function $f(t)$ of the integrand can be expanded in powers of $\zeta$. This gives the computational scheme for the coefficients $c_k$ when $\vert\zeta\vert$ is small. For other values of $\zeta$ the removable singularities in $c_k$ are not of any numerical concern.

The transformations of the integrals into Airy-type integrals considered in the previous section are principal tools for obtaining uniform Airy-type expansions. Because these asymptotic techniques are not the starting point of the present paper, in this section  we will consider the numerical evaluation of an integral of the $J$-Bessel  function instead of applying the cubic transformation \eqref{eq:transBes03}.

We return to the integral in  \eqref{eq:transBes02}. When $0<z<1$, the saddle point contour through the positive saddle point $s_+=\arccosh (1/z)$  is given by (we write $s=\sigma+i\tau$)
\begin{equation}\label{eq:Jbesor01}
z\cosh\sigma \sin\tau-\tau=0\quad \Longrightarrow \quad \sigma=\arccosh \frac{\tau}{z\sin\tau},\quad
-\pi<\tau<\pi.
\end{equation}
When we use this parametrization for the saddle point contour, we can write the integral representation in \eqref{eq:transBes02} in the form
\begin{equation}\label{eq:Jbesor02}
J_\nu(\nu z)=\frac{1}{2\pi i}\int_{-\pi}^\pi e^{\nu (z\sinh\sigma\cos\tau-\sigma)}\left(\frac{d\sigma}{d\tau}+i\right)\,d\tau.
\end{equation}
Because $\sigma$ is an even function of $\tau$, we have
\begin{equation}\label{eq:Jbesor03}
J_\nu(\nu z)=\frac{e^{-\nu\rho}}{2\pi }\int_{-\pi}^\pi e^{-\nu\psi(\tau)} \,d\tau,
\end{equation}
where, with $\sigma$ as defined in \eqref{eq:Jbesor01},
\begin{equation}\label{eq:Jbesor04}
\begin{array}{@{}r@{\;}c@{\;}l@{}}
\psi(\tau)&=&-\rho-z\sinh\sigma\cos\tau+\sigma, \\[8pt]
\rho&=&\arccosh(1/z)-\sqrt{1-z^2}= \arctanh\sqrt{1-z^2}-\sqrt{1-z^2}.
\end{array}
\end{equation}
An expansion for small values of $\vert\tau\vert$ reads
\begin{equation}\label{eq:Jbesor05}
\psi(\tau)=\tfrac{1}{2}\sqrt{1-z^2}\,\tau^2-\frac{3z^2+2}{72 \sqrt{1-z^2}}\tau^4+\bigO\left(\tau^5\right).
\end{equation}
This is valid for $0<z<1$. When $z=1$ we have 
\begin{equation}\label{eq:Jbesor06}
\psi(\tau)=\tfrac{4}{27}\sqrt{3}\,\vert\tau\vert^3+\tfrac{8}{14175} \sqrt{3}\,\vert\tau\vert^7+\bigO\left(\tau^9\right).
\end{equation}

\begin{figure}
\begin{center}
\epsfxsize=10cm \epsfbox{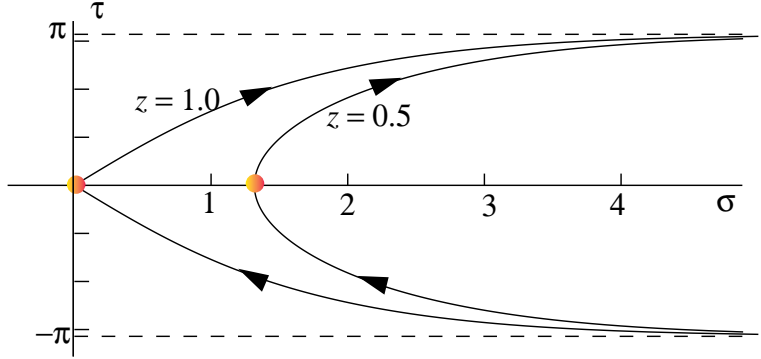}
\caption{
\label{fig:fig06} Saddle point contours defined by \eqref{eq:Jbesor01} in the $s$-plane, $s=\sigma+i\tau$ for $z=0.5$ and $z=1.0$.  }
\end{center}
\end{figure}

In Figure~\ref{fig:fig06} we show the saddle point contours defined by \eqref{eq:Jbesor01} in the $s$-plane, $s=\sigma+i\tau$ for $z=0.5$ and $z=1.0$. When $z=0.5$ there is a smooth passage through the saddle point at $s_+\doteq 1.32$, when $z=1.0$ the contour shows a kink at the saddle point $s_+=0$.

We can evaluate the integral in \eqref{eq:Jbesor03} using the trapezoidal rule on a finite interval, as long as $z\in(0,1)$ is not near 1. We see from the expansions of $\psi(\tau)$ that this function is not analytic at $\tau=0$ when $z\to1$.  

When $z>1$ we can use the Hankel function representation in the form (see \eqref{eq:etaneg13})
\begin{equation}\label{eq:Jbesor07}
H_\nu^{(1)}(\nu z)=\frac{1}{\pi i}\int_{-\infty}^{\infty +\pi i}e^{\nu (z\sinh s-s)}\,ds,
\end{equation}
and we use the contour through the saddle point $is_+$, with $s_+=\arccos(1/z)$. Again we write $s=\sigma+i\tau$, and the contour of steepest descent through $s=is_+$ is defined by 
\begin{equation}\label{eq:Jbesor08}
\Im(z\sinh s-s)=\Im(\sinh(is_+)-is_+)=\wt\rho,
\end{equation}
or
\begin{equation}\label{eq:Jbesor09}
z\cosh\sigma\sin\tau-\tau=\wt\rho,
\end{equation}
where
\begin{equation}\label{eq:Jbesor10}
\begin{array}{@{}r@{\;}c@{\;}l@{}}
\wt\rho&=&z\sin(s_+)-s_+\\[8pt]
&=&\sqrt{z^2-1}-\arccos(1/z)=\sqrt{z^2-1}-\arctan\sqrt{z^2-1}.
\end{array}
\end{equation}

\begin{figure}
\begin{center}
\epsfxsize=10cm \epsfbox{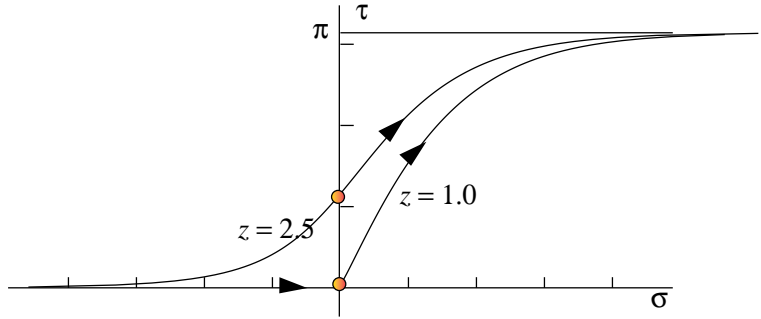}
\caption{
\label{fig:fig07} Saddle point contours in the $s$-plane defined by \eqref{eq:Jbesor11}, $s=\sigma+i\tau$, for $z=1.0$ and $z=2.5$. }
\end{center}
\end{figure}

From \eqref{eq:Jbesor09} we find for the contour the representation
\begin{equation}\label{eq:Jbesor11}
 \sigma=\arccosh \frac{\tau+\wt\rho}{z\sin\tau},\quad
0<\tau<\pi, \quad \sign(\sigma)=\sign(\tau-s_+),
\end{equation}
and for the integral
\begin{equation}\label{eq:Jbesor12}
H_\nu^{(1)}(\nu z)=\frac{e^{i\nu\wt\rho}}{\pi i}\int_{0}^{\pi }e^{-\nu \wt\psi(\tau)}\left(\frac{d\sigma}{d\tau}+i\right)\,d\tau,
\end{equation}
where
\begin{equation}\label{eq:Jbesor13}
 \wt\psi(\tau)=-z\sinh\sigma\cos\tau+\sigma,\quad 
 \frac{d\sigma}{d\tau}=\frac{1-z\cosh\sigma\cos\tau}{z\sinh\sigma\sin\tau}.
 \end{equation}
For small values of $\vert\tau-s_+\vert$ we have the expansion
\begin{equation}\label{eq:Jbesor14}
\wt\psi(\tau)=\sqrt{z^2-1}\,(\tau-s_+)^2+\frac{3z^2+2}{18 \sqrt{z^2-1}}(\tau-s_+)^4+\bigO\left((\tau-s_+)^5\right).
\end{equation}
This is valid when $z>1$. 

The integrand of the integral in \eqref{eq:Jbesor12} has its maximal value at $\tau=s_+=\arccos(1/z)$, and for this value $\sigma=0$. The integrand vanishes at the endpoints with all its derivatives, and, again, we can use the trapezoidal rule for numerical evaluations, when $z>1$, but $z$ should not be close to 1.

When $\nu$ and $z$ are real, the result for $J_\nu(\nu z)$ follows by taking the real part of the computed value of $H_\nu^{(1)}(\nu z)$.

In Figure~\ref{fig:fig07} we show the saddle point contours defined by \eqref{eq:Jbesor11} in the $s$-plane, $s=\sigma+i\tau$ for $z=1.0$ and $z=2.5$. When $z=2.5$ the passage through the saddle point at $s_+\doteq 1.16 i$ runs smoothly. When $z=1.0$ the negative axis is part of the contour, and  the contour shows a kink at the saddle point $s_+=0$.

Finally we consider the case $z\sim1$. The saddle point contours for the cases $0<z<1$ (see Figure~\ref{fig:fig06})  and for $z > 1$  (see Figure~\ref{fig:fig07})  become non-smooth when $z\to 1$. In that case $\rho$ defined in \eqref{eq:Jbesor04} and $\wt\rho$ defined in \eqref{eq:Jbesor10} become small, and the corresponding
exponential factors in  front of the integral may become less important.  Of course, this depends on the value of $\nu$.

When $\nu\rho$ is not large, say, $\nu\rho\le 1$, it is not needed to use the saddle point contour, nor the saddle point. We can use, for example, the contour depicted in Figure~\ref{fig:fig06} for $z=0.5$, and shift it through the saddle point $s_+$. In this way, when $0<z\le1$, we replace the contour defined in \eqref{eq:Jbesor01} by the contour defined by 
\begin{equation}\label{eq:Jbesor15}
\sigma=\arccosh\frac{2\tau}{\sin \tau}-\arccosh(2) + s_+,\quad -\pi<\tau<\pi.
 \end{equation}
For small values of $\tau$ we have the expansion
\begin{equation}\label{eq:Jbesor16}
\sigma=s_++\tfrac19\sqrt{3}\,\tau^2\left(1+\tfrac1{180}\tau^2+\tfrac{37}{7560}\tau^4+\ldots\right).
 \end{equation}

The representation in \eqref{eq:transBes02} becomes
\begin{equation}\label{eq:Jbesor17}
J_\nu(\nu z)=\frac{1}{2\pi}\int_{-\pi}^{\pi}e^{ -\nu p(\tau)} q(\tau)\,d\tau,
\end{equation}
where 
\begin{equation}\label{eq:Jbesor18}
\begin{array}{@{}r@{\;}c@{\;}l@{}}
p(\tau)&=&\sigma-z\sinh\sigma\cos\tau, \\[8pt]
q(\tau)&=&\dsp{ \cos r(\tau)+\sin r(\tau)\,\frac{d\sigma}{d\tau},}\\[8pt]
r(\tau)&=&\nu\left(z\cosh\sigma\sin\tau-\tau\right), \\[8pt]
\dsp{\frac{d\sigma}{d\tau}}&=&2\dsp{\frac{\sin\tau-\tau\cos\tau}{\sin^2\tau \sinh\left(\frac{2\tau}{\sin \tau}\right)}.} 
\end{array}
\end{equation}

For $z\ge1$, as long as $-1\le \eta\le0$, we modify the contour in \eqref{eq:Jbesor15} by writing
\begin{equation}\label{eq:Jbesor19}
\sigma=\arccosh\frac{2\tau}{\sin \tau}-\arccosh\frac{2z\arccos(1/z)}{\sqrt{z^2-1}},\quad -\pi<\tau<\pi.
\end{equation}
With this choice of the parameter $\sigma$  the contour runs through the two points $(\sigma,\tau)=(0,\pm \arccos(1/z))$, that is, through the two saddle points on the imaginary axis. The formulas \eqref{eq:Jbesor17}--\eqref{eq:Jbesor18} remain the same.
 
\begin{remark}\label{rem:rem05}
In the Airy-type integrals we have decided between the methods by verifying the value of $\eta$. To use this for the integrals for the Bessel function in the present section, we observe that when $0<z<1$ the relation between $\rho$ in \eqref{eq:Jbesor04} and $\zeta$ in 
\eqref{eq:transBes04} reads $\rho=\frac23\zeta^{\frac32}$ and the corresponding $\eta$ in \eqref{eq:transBes04} satisfies $\eta=\nu^{\frac23}\zeta$, that is, $\eta=\left(\frac32\nu\rho\right)^{\frac23}$. For $z>1$ we can replace $\rho$ by  $\wt\rho$ given in \eqref{eq:Jbesor10}.
\eoremark
\end{remark}

\section{Numerical examples}\label{sec:numerics}
Earlier in the text we have given a few examples of the numerical performance of the trapezoidal rule. We give results for the 
integrals derived in \S\ref{sec:airytype}  and for the integral representations of the Bessel function in \S\ref{sec:org}. The computations are done by using Maple with Digits=16. For comparison we used the codes for the Airy function and $J$-Bessel function of Maple, and used for that purpose Digits =24.

\begin{table}
\caption{
Relative errors $\delta$ in the computation of  the integral \eqref{eq:intro01} with  $f(t)=\cos t$ for several values of $\eta$.  For more details, see the text.
\label{tab:tabl02}}
$$
\begin{array}{rcc | rcc | ccc}
\eta &k_m&\delta&\eta \ \ &k_m&\delta&\eta &k_m&\delta\  \\
\hline
-1  &  49  &  0.18{\rm e-}14      &   -1.00    &  32  &   0.35{\rm e-}14     & 1  &  37  &  0.10{\rm e-}14  \\
-2  &  38  &   0.18{\rm e-}14     &   -0.60      &  32  &   0.30{\rm e-}14     & 2  &  22  &  0.10{\rm e-}14 \\
-3  &  32  &   0.28{\rm e-}12     &   -0.20     &  31  &   0.26{\rm e-}14     & 3  &  16  &  0.00{\rm e-}00\\
-4  &  28  &   0.72{\rm e-}10     &    0.20        &  31  &   0.40{\rm e-}15     & 4  &  13  &  0.30{\rm e-}13\\
-5  &  25  &  0.63{\rm e-}10      &    0.60      &  31  &  0.00{\rm e-}00      &  5  &  11  &  0.42{\rm e-}11\\
-6  &  22  &  0.16{\rm e-}07      &   1.00        &  31  &   0.18{\rm e-}14    &   6  &  10  &  0.28{\rm e-}09\\
\hline
\end{array}
$$
\end{table}

\subsection{Integrals derived in \S\ref{sec:airytype}}\label{sec:intairy}

In Table~\ref{tab:tabl02} we give the relative errors $\delta$ in the computation of  the integral \eqref{eq:intro01} for several values of $\eta$. We take $f(t)=\cos t$ and we give the number of terms $k_m$ needed to have the integrand values less than $1.0{\rm e-}16$ for $k=k_m$. Because the exponential function in the integrands in \eqref{eq:etalarge07}, \eqref{eq:etaneg11} and \eqref{eq:etasmall03} are quite different with respect to fast convergence, we take different stepsizes for each interval. We take $h=0.2$ if $\eta \le -1$, $h=0.05$ if $-1\le \eta \le 1$, and $h=0.3$ if $\eta \ge 1$. 

We see from Table~\ref{tab:tabl02} that a fixed stepsize $h$ for each interval is not a guarantee of good performance for all $\eta$, except in the middle interval.

\subsection{Integral representations of the Bessel function in \S\ref{sec:org}}\label{sec:intbess}

We have tested the trapezoidal rule for several cases. 
We computed the errors by computing three successive Bessel functions and  verified the recurrence relation
(see \cite[\S10.6]{Olver:2010:BFS})
\begin{equation}\label{eq:numer01}
2\nu J_\nu(z)=z\left(J_{\nu-1}(z)+J_{\nu+1}(z)\right).
\end{equation}
In Table~\ref{tab:tabl03} and Table~\ref{tab:tabl04} we also compare the answer with the computation by the Maple code for $J_\nu(x)$, with Digits $=24$. These errors are in the column $J_\nu$-error. For the extreme high $\nu$-values in  Table~\ref{tab:tabl04} we only verified by recursion.
\begin{enumerate} 
\item
We used the method described in \S\ref{sec:etasmall}, and verified if we could use it not only for $-1\le \eta \le 1$ but also for complex values of $\eta$ inside the unite circle. Indeed, the method can be used without further preparations.  In a numerical example, we have computed the Airy function $\Ai(\eta)$  by taking in \eqref{eq:etasmall02} $f(t)=1$, and we used the trapezoidal rule with $h=0.06$ for $\eta=e^{k\pi i/16}$ with $k=0,1,2,\ldots,16$. We found a maximal absolute value $6.24\times10^{-15}$ of the absolute  error for $k=14$. Computations were done in Maple with Digits $=16$. We summed the series of the trapezoidal rule for the integral in \eqref{eq:etasmall02} until the absolute value of the terms divided by the sum obtained so far became less than $10^{-16}$.  In this way, and using that the integrand is an even function of $\theta$, 26 terms were needed.

\begin{table}
\caption{
Relative errors in the computation of $J_\nu(x)$  for $\nu=100$ and $x=91, 93,\ldots,99$, by using the integral representation  in \eqref{eq:Jbesor03}. The method is intended for $\eta\ge1$, and we see bad performance for the final $x$-value because one of the  functions in the recursion relation has $\nu=x=99$, which for that function $\eta=0$.
\label{tab:tabl03}}
$$
\begin{array}{rcccc}
x & \eta &J_\nu(x)&\   {\rm  rec. \ error}\ & \  J_\nu-\rm{error}\   \\
\hline
91 & 2.51 & 0.4256251712037803{\rm e-}2 & 0.20{\rm e-}14 & 0.22{\rm e-}14\\
93 & 1.94 & 0.1050032579531836{\rm e-}1 & 0.23{\rm e-}14 & 0.18{\rm e-}14 \\
95 & 1.38 & 0.2315076800942791{\rm e-}1 & 0.43{\rm e-}14 & 0.24{\rm e-}14\\
97 & 0.82 & 0.4528109693556812{\rm e-}1 & 0.60{\rm e-}14 & 0.16{\rm e-}14 \\
99 & 0.27 & 0.7768716170045931{\rm e-}1 & 0.27{\rm e-}05 & 0.12{\rm e-}14\\
\hline
\end{array}
$$
\end{table}

\item
We have computed $J_\nu(x)$ for $\nu=100$ and $x=90, 91,\ldots, 99$, by using the integral representation  in \eqref{eq:Jbesor03}.  The results are shown in Table~\ref{tab:tabl03}. The stepsize for the trapezoidal rule is $h=0.05$, and the number of terms is 24. The corresponding $\eta$ values (see Remark~\ref{rem:rem05}) are not always larger than 1, and for $x=99$,  we see bad performance for the final $x$-value because one of the  functions in the recursion relation has $\nu=x=99$, and for that function $\eta=0$.

\begin{table}
\caption{
Relative errors in the computation of $J_\nu(x)$  for $\nu=100$ and $x=99.0, 99.2,\ldots,100$, by using the integral representation  in \eqref{eq:Jbesor17}. 
\label{tab:tabl04}}
$$
\begin{array}{rcccc}
x & \eta &J_\nu(x)&\   {\rm  rec. \ error}\ & \  J_\nu-\rm{error}\   \\
\hline
99.0 & 0.272 & 0.7768716170045941{\rm e-}1 & 0.18{\rm e-}14 & 0.12{\rm e-}15 \\
99.2 & 0.215 & 0.8135695322732582{\rm e-}1 & 0.73{\rm e-}15 & 0.16{\rm e-}14\\
99.4 & 0.163 & 0.8507190689984157{\rm e-}1 & 0.56{\rm e-}15 & 0.28{\rm e-}14\\
99.6 & 0.109 & 0.8882046195955568{\rm e-}1 & 0.40{\rm e-}14 & 0.17{\rm e-}14\\
99.8 & 0.054 & 0.9258996685877174{\rm e-}1 & 0.18{\rm e-}14 & 0.40{\rm e-}14\\
100.0  & 0.000 & 0.9636667329586151{\rm e-}1& 0.23{\rm e-}15 & 0.52{\rm e-}15\\
\hline
\end{array}
$$
\end{table}

\item
 In Table~\ref{tab:tabl04}  we show the results for $0<\eta<1$, with $\nu=100$ and  $x=99.1, 99.2,\ldots,100$. We used  the integral representation  in \eqref{eq:Jbesor17}.  The stepsize for the trapezoidal rule is $h=0.05$, and the number of terms is~28.

\begin{table}
\caption{
Relative errors in the computation of $J_\nu(\nu z)$ by using the trapezoidal for  the integral in \eqref{eq:transBes08} for rather extreme values of $\nu$ ($\nu=10^k$). More details are given in the text.
\label{tab:tabl05}}
$$
\begin{array}{rccccc}
k & z & \zeta & J_\nu(\nu z)&  {\rm rec.\   error} & J_\nu-\rm{error}\\
\hline
2  &  0.927948934   &  0.93{\rm e-}1   &  0.9620266889434034{\rm e-}2   &  0.69{\rm e-}15   & 0.40{\rm e-}14 \\
4  &  0.996583557   &  0.43{\rm e-}2   &  0.2043772855795365{\rm e-}2   &  0.21{\rm e-}14   & 0.00{\rm e-}00\\
6  &  0.999841268   &  0.20{\rm e-}3   &  0.4400304405124362{\rm e-}3   &  0.74{\rm e-}14   &0.50{\rm e-}14 \\
8  &  0.999992632  &   0.93{\rm e-}5   &  0.9479881456179256{\rm e-}4   &  0.86{\rm e-}14   &0.33{\rm e-}14 \\
10&  0.999999658   &  0.43{\rm e-}6   &  0.2042375676682798{\rm e-}4   &  0.42{\rm e-}14   & 0.20{\rm e-}14\\
\hline
\end{array}
$$
\end{table}

\item

For Table~\ref{tab:tabl05} we have used the values $\eta=2$, $\nu=10^k$, $k=2,4\ldots,10$ and 16 terms (with positive index $k$) in the 
trapezoidal rule for the Airy-type integral representation of $J_\nu(\nu z)$  in \eqref{eq:transBes08}. 
 We give values of the corresponding $z$, $\zeta$, $J_\nu(\nu z)$, and the relative errors based on the recursion and on an algorithm  described in \cite{Temme:1997:NAU}.
We show less than 16 relevant digits for $z$ and $\zeta$  to keep the table in a proper size. When $\nu$ and $\eta$ are given, $\zeta$ follows from \eqref{eq:transBes08}, that is, $\zeta=\eta\nu^{-\frac{2}{3}}$, and $z$ follows from the first line in \eqref{eq:transBes04}. We observe for large values of $\nu$  the small values of $\zeta$ and values of $z$ close to unity when we use $\eta=2$. In this numerical example for large $\nu$, with $\eta$ and $\nu$ given, it is important first to find $1-z$, which is of order $\zeta$, and then $z$ (see also expansion \eqref{eq:transBes04}). In several formulas the accuracy of $1-z$ is relevant. For example, when computing the saddle point $s_+=\arccosh(1/z)$, which can also be written as 
$s_+=\arctanh\left(\sqrt{1-z^2}\right)$.

\end{enumerate}


\section*{Acknowledgments}
We thank the referees for their constructive and helpful remarks. \\
We acknowledge financial support from Ministerio de Ciencia e Innovaci\'on, Spain, 
projects MTM2015-67142-P (MINECO/FEDER, UE) and PGC2018-098279-B-I00 (MCIU/AEI/FEDER, UE). \\
NMT thanks CWI, Amsterdam, for scientific support.

\bibliographystyle{plain}

\begin{thebibliography}{10}

\bibitem{Berry:2010:IWC}
M.~V. Berry and C.~J. Howls.
\newblock Chapter 36, {I}ntegrals with coalescing saddles.
\newblock In {\em N{IST} {H}andbook of {M}athematical {F}unctions}, pages 775--793.
  U.S. Dept. Commerce, Washington, DC, 2010.
\newblock http://dlmf.nist.gov/36.

\bibitem{Borghi:2016:COS}
R.~Borghi.
\newblock Computational optics through sequence transformations.
\newblock In {\em Progress in {O}ptics, Volume 61}, pages 1--68. Elsevier B.V.,
  Amstrdam, 2016.

\bibitem{Chester:1957:ESD}
C.~Chester, B.~Friedman, and F.~Ursell.
\newblock An extension of the method of steepest descents.
\newblock {\em Proc. Cambridge Philos. Soc.}, 53:599--611, 1957.

\bibitem{Connor:1990:PMU}
J.~N.~L. Connor.
\newblock Practical methods for the uniform asymptotic evaluation of
  oscillating integrals with several coalescing saddle points.
\newblock In {\em Asymptotic and computational analysis ({W}innipeg, {MB},
  1989)}, volume 124 of {\em Lecture Notes in Pure and Appl. Math.}, pages
  137--173. Dekker, New York, 1990.

\bibitem{Connor:1992:UAO}
J.~N.~L. Connor, P.~R. Curtis, and R.~A.~W. Young.
\newblock Uniform asymptotics of oscillating integrals: applications in
  chemical physics.
\newblock In {\em Wave asymptotics ({M}anchester, 1990)}, pages 24--42.
  Cambridge Univ. Press, Cambridge, 1992.

\bibitem{Deano:2009:CGQ}
A.~Dea\~{n}o and D.~Huybrechs.
\newblock Complex {G}aussian quadrature of oscillatory integrals.
\newblock {\em Numer. Math.}, 112(2):197--219, 2009.

\bibitem{Deano:2018:CHO}
A.~Dea\~{n}o, D.~Huybrechs, and A.~Iserles.
\newblock {\em Computing highly oscillatory integrals}.
\newblock Society for Industrial and Applied Mathematics (SIAM), Philadelphia,
  PA, 2018.

\bibitem{Dunster:2017:CAE}
T.~M. Dunster, A.~Gil, and J.~Segura.
\newblock Computation of asymptotic expansions of turning point problems via
  {C}auchy's integral formula: {B}essel functions.
\newblock {\em Constr. Approx.}, 46(3):645--675, 2017.

\bibitem{Ferreira:2018:AES}
C.~Ferreira, J.~L. L\'{o}pez, and E.~P\'erez-Sinusia.
\newblock The asymptotic expansion of the swallowtail integral in the highly
  oscillatory region.
\newblock {\em Appl. Math. Comput.}, 339:837--845, 2018.

\bibitem{Gil:2001:ONI}
A.~Gil, J.~Segura, and N.~M. Temme.
\newblock On nonoscillating integrals for computing inhomogeneous {A}iry
  functions.
\newblock {\em Math. Comp.}, 70(235):1183--1194, 2001.

\bibitem{Gil:2002:AGH}
A.~Gil, J.~Segura, and N.~M. Temme.
\newblock Algorithm 822: {GIZ}, {HIZ}: two {F}ortran 77 routines for the
  computation of complex {S}corer functions.
\newblock {\em ACM Trans. Math. Software}, 28(4):436--447, 2002.

\bibitem{Gil:2007:NSF}
A.~Gil, J.~Segura, and N.~M. Temme.
\newblock {\em Numerical methods for special functions}.
\newblock Society for Industrial and Applied Mathematics (SIAM), Philadelphia,
  PA, 2007.

\bibitem{Gil:2002:CCA}
Amparo Gil, Javier Segura, and Nico~M. Temme.
\newblock Computing complex {A}iry functions by numerical quadrature.
\newblock {\em Numer. Algorithms}, 30(1):11--23, 2002.

\bibitem{Goodwin:1949:EIF}
E.~T. Goodwin.
\newblock The evaluation of integrals of the form $\int\sp \infty\sb {-\infty}
  f(x) e\sp {-x\sp {2}} dx$.
\newblock {\em Proc. Cambridge Philos. Soc.}, 45(2):241--245, 1949.

\bibitem{Huybrechs:2018:ANM}
D.~Huybrechs, A.~B.~J. Kuijlaars, and N.~Lejon.
\newblock A numerical method for oscillatory integrals with coalescing saddle
  points.
\newblock {\em SIAM J. Numer. Anal.}, 57(6):2707--2729, 2019.

\bibitem{Kirk:2000:NEO}
N.~P. Kirk, J.~N.~L. Connor, and C.~A. Hobbs.
\newblock An adaptive contour code for the numerical evaluation of the
  oscillatory cuspoid canonical integrals and their derivatives.
\newblock {\em Comput. Phys. Commun.}, 132:142--165, 2000.

\bibitem{Lopez:2017:AFE}
J.~L. L\'{o}pez and P.~J. Pagola.
\newblock Analytic formulas for the evaluation of the {P}earcey integral.
\newblock {\em Math. Comp.}, 86(307):2399--2407, 2017.

\bibitem{Milovanovic:2017:CIH}
G.~V. Milovanovi\'{c}.
\newblock Computing integrals of highly oscillating special functions using
  complex integration methods and {G}aussian quadratures.
\newblock {\em Dolomites Res. Notes Approx.}, 10(Special Issue):79--96, 2017.

\bibitem{Olver:2010:ARF}
F.~W.~J. Olver.
\newblock Chapter 9, {A}iry and related functions.
\newblock In {\em N{IST} {H}andbook of {M}athematical {F}unctions}, pages
  193--213. U.S. Dept. Commerce, Washington, DC, 2010.
\newblock http://dlmf.nist.gov/9.

\bibitem{Olver:2010:BFS}
F.~W.~J. Olver and L.~C. Maximon.
\newblock Chapter 10, {B}essel functions.
\newblock In {\em N{IST} {H}andbook of {M}athematical {F}unctions}, pages
  215--286. U.S. Dept. Commerce, Washington, DC, 2010.
\newblock http://dlmf.nist.gov/10.

\bibitem{Temme:1997:NAU}
N.~M. Temme.
\newblock Numerical algorithms for uniform {A}iry-type asymptotic expansions.
\newblock {\em Numer. Algorithms}, 15(2):207--225, 1997.

\bibitem{Temme:2015:AMI}
N.~M. Temme.
\newblock {\em Asymptotic methods for integrals}, volume~6 of {\em Series in
  Analysis}.
\newblock World Scientific Publishing Co. Pte. Ltd., Hackensack, NJ, 2015.

\bibitem{Trefethen:2014:TEC}
L.~N. Trefethen and J.~A.~C. Weideman.
\newblock The exponentially convergent trapezoidal rule.
\newblock {\em SIAM Rev.}, 56(3):385--458, 2014.

\end{thebibliography}

\end{document}